\documentclass[10pt]{amsart}
\usepackage{a4,amsmath,amssymb,amsthm}
\usepackage[english]{babel}

\newtheorem{theorem}{Theorem}[section]
\newtheorem{definition}[theorem]{Definition}

\newtheorem{lemma}[theorem]{Lemma}
\newtheorem{proposition}[theorem]{Proposition}
\newtheorem{corollary}[theorem]{Corollary}

\theoremstyle{remark}
\newtheorem{remark}[theorem]{Remark}
\newtheorem{example}[theorem]{Example}
\numberwithin{equation}{section}

\begin{document}

\title{INVARIANT MEASURES FOR INTERVAL TRANSLATIONS AND SOME OTHER PIECEWISE CONTINUOUS MAPS}
\author{Sergey Kryzhevich}

\begin{abstract} We study some special classes of piecewise continuous maps on a finite smooth partition of a compact manifold and look for invariant measures for such maps.  
We show that in the simplest one-dimensional case (so-called interval translation maps) a Borel probability non-atomic invariant measure exists for any map.  We use this result to demonstrate that any interval translation map endowed with such a measure is metrically equivalent to an interval exchange map. Finally, we study the general case of piecewise continuous maps and prove a simple result on existence of an invariant measure provided all discontinuity points are wandering. \end{abstract}
\maketitle

\section*{Introduction}

The main objective of this research is to understand how the numerical methods could be used to study piecewise continuous dynamics. Modelling the discontinuous dynamics manifests many difficulties that could be easily illustrated by the following example \cite[Exercise 4.1.1]{KH}.

The famous Krylov - Bogolyubov Theorem claims that any continuous transformation of a compact metric space admits a Borel probability invariant measure. The similar statement for discontinuous maps is wrong.

\begin{example} Consider the map $T:[0,1]\to [0,1]$ given by the formula: $T(x)=x/2$ if $x>0$; $T(0)=1$. This map does not admit any Borel probability invariant measure.
\end{example}

So, the numerical methods may be inappropriate even in their weakest form -- modelling invariant measures. Observe that for continuous maps of compact sets we have lower semicontinuity of the set of invariant measures that fails for discontinuous maps. 

Applying numerical methods, one usually cuts the phase space into small indivisible pieces (pixels, finite elements etc) that do not change their shapes. The numerical method shuffles these pieces, sometimes with overlaps. Here we face at least two questions:
\begin{itemize}
\item[Q1:] do such piecewise isometries admit any invariant measures?
\item[Q2:] do these measures approximate any invariant measure of the initial system? 
\end{itemize}

The first problem is well-known in the one-dimensional case.  

\begin{example} Consider the circle ${\mathbb T}^1:={\mathbb R}/{\mathbb Z}$. We represent it as a union of disjoint subsegments $M_{j}=[t_j,t_{j+1})$, $j=0,\ldots,n$, $t_0=t_n$ and define the map $S$ by the formula 
$$S(t)=t+c_j \mod 1, \qquad t\in M_{j}.$$
Here $c_j$ are real values. Such map is called \emph{interval translation} (ITM) or, if it is one-to-one it is called \emph{interval exchange} (IEM).   
\end{example}

Similarly one may consider interval translation maps on the segment $[0,1]$. Notice that in this paper we consider the orientation - preserving maps only.  

For interval exchange maps the question Q1 is trivial: the Lebesgue measure is always invariant. Moreover, the map $S$ admits at most $n$ Borel probability invariant non-atomic ergodic measures (see \cite[\S 14.5, \S14.6]{KH} for the basic theory of interval exchange maps and \cite{V1} and \cite{Y} as surveys for deeper results).

The case of non-invertible ITMs, first considered by M.\, Boshernitzan and I.\, Kornfeld \cite{BK} is much more sophisticated. One of the principal problems for ITMs is their classification. 

\begin{definition}
 We say that the ITM $S$ is \emph{finite} if there exists a number $m\in {\mathbb N}$ such that 
 $$S^m({\mathbb T}^1)= \bigcap_{k=1}^\infty S^k ({\mathbb T}^1).$$
 Otherwise, the map $S$ is \emph{infinite}. 
\end{definition}

In \cite{BK}, authors have demonstrated that many ITMs are finite and thus may be restricted to interval exchange maps. However, there are examples with ergodic measures  supported on Cantor sets. J.\,Schmeling and S.Troubetzkoy \cite{ST} provided some estimates on the number of minimal subsets for ITMs. 

H.\,Bruin and S.\, Troubetzkoy \cite{BT} studied ITMs of a segment of 3 intervals ($n=3$). It was shown that in this case typical ITM is finite. In any case, results on Hausdorff dimension for attractors and unique ergodicity are given. These results are generalised in \cite{B2} for ITMs with arbitrary many pieces. There is an uncountable set of parameters leading to type $\infty$ interval translation maps but the Lebesgue measure of these parameters is zero. Furthermore conditions are given that imply that the ITMs have multiple ergodic invariant measures. H. Bruin and G. Clark \cite{BC} studied the so-called double rotations ($n=2$ for maps of the circle). Almost all double rotations are of finite type. The parameters that correspond to infinite type maps, form a set of Hausdorff dimension strictly between 2 and 3. 

J.\, Buzzi and P.\, Hubert \cite{BH} studied piecewise monotonous maps of zero entropy and no periodic points. Particularly, they demonstrated that orientation-preserving ITMs without periodic points may have at most $n$ ergodic probability invariant measures where $n$ is the number of intervals.  

D.\,Volk \cite{V2} was studying ITMs of the segment. He demonstrated that almost any (w.r.t. Lebesgue measure on the parameters set) ITM of 3 intervals is conjugated to a rotation or to a double rotation and, hence, is of finite type.

B.\, Pires in his preprint \cite{P} proved that almost all ITMs admit a non-atomic invariant measure (he assumed that the map does not have any connections or periodic points). In our paper, we generalise the referred result proving the same result for all ITMs. Moreover, our techniques allow us to claim that all restrictions of ITMs to supports of their invariant measures are metrically equivalent to interval exchange maps. A very similar statement was proved in \cite{P}: any injective piecewise continuous map $S:[0,1]\to [0,1]$ is semi-conjugate to an interval exchange transformation, possibly with flips.

The more general class of maps is the so-called piecewise translation maps (PWT), see the precise definition below. This is the multi-dimensional generalisation of ITMs. Such maps are widely used in applications: herding dynamics in Markov networks, second order digital filters, sigma-delta modulators, buck converters, three-capacitance models, error diffusion algorithm in digital printing, machine learning etc, please see \cite{SIA}, \cite{T} and \cite{V2} for review and references. 

A good survey on the results on the theory of PWTs and some brilliant original results were given by D.\, Volk \cite{V3}. Particularly, such properties as ergodicity, finite type of the map and the structure of the attractor were discussed. A very interesting result was obtained in the preprint \cite{T}. If the dimension of the phase space $d$ and the number $n$ of the elements $M_j$ of the partition are such that $n=d+1$, it was proved that any PWT with rationally independent shift vectors is of finite type and, consequently, the Lebesgue measure of the attractor is non-zero and, hence there exist a continuous invariant measure. This result was proved by the powerful techniques of restriction of the discontinuous multivalued map to the map on the set of compact subsets. A similar approach was used in the preprint \cite{V3}.

The most general class of the maps is the so-called piecewise isometries. We split the phase space into a finite number of parts and assume that all restrictions of the considered maps are isometries (that may include rotations).
 
J.\,Buzzi \cite{B1} demonstrated that piecewise isometries defined on a finite union of polytopes have zero topological entropy in any dimension. Xin-chu Fu and Jinqiao Duan \cite{FD} studied this problem in dimension 2 (the so-called planar isometries). They provided sufficient conditions for existence of Milnor-type attractors. Zhan-he Chen and Rong-zhong Yu \cite{CYF} have also considered planar piecewise isometries. 
 
The partition of the phase space engenders a naturally defined symbolic dynamics. A. Goetz \cite{G1,G2} demonstrated that for so-called regular partition (that is true under assumptions of our paper) symbolic dynamics of the isometry cannot embed subshifts of finite type with positive entropy. The condition of polynomial growth of symbolic words is given.

In this paper we do the following. We formulate conditions on the piecewise continuous maps that are sufficient for existence of Borel probability invariant measures. Particularly, we demonstrate that any ITM admits a non-atomic Borel probability invariant measure and is metrically equivalent to an interval exchange map. Finally, we study the case when invariant measures of piecewise continuous maps may be approximated by ones of approximating PWTs. 

Later on we use designations $f(t_0-0)$ and $f(t_0+0)$ for the left and the right limit respectively.

\section{Piecewise isometries, their periodic points and domains.}

Consider a $d$-dimensional compact Riemannian manifold $(M,{\mathrm{dist}\,})$ or a bounded set $M\subset {\mathbb R}^d$ that can be represented as a finite union of closures of pairwise disjoint domains $M_j$ ($j=1,\ldots, n$). We consider a map $S:M\to M$ such that all restrictions $S|_{\mathrm {Int}\, M_k}$ are isometries (we call such map a \emph{piecewise isometry}) or translations (the map is called \emph{piecewise translation}).

Let $H=\bigcup_{j=1}^n {\partial M_j}$. Using terminology of Piecewise Continuous Dynamics, we call $H$ \emph{discontinuity set} for the map $S$, let $\Omega:=\bigcup_{k \ge 0} S^{-k} (H)$ be the set of all eventually discontinuity points.

On $M$, we consider the Lebesgue probability measure $\mathrm{Leb}$. We assume that $\mathrm{Leb}\,(H)=0$ that implies $\mathrm{Leb}\,(\Omega)=0$.  

Let us study invariant measures and invariant sets of the map $S$. We use ideas of \cite[\S14.5]{KH}. There, properties of the so-called interval exchange maps are discussed. In fact, "regular" interval exchange maps are a particular case of maps, we are studying. However, our case is much more difficult and many techniques of the IEM theory cannot be applied directly.  

Observe that the measure $\mathrm{Leb}$ is not invariant for $S$ unless this map is invertible almost everywhere. Moreover, generally speaking, it is not evident if all interval translation maps admit Borel probability invariant measures. The positive answer to this question is one of the main results of the paper. 

\begin{definition} We say that a subset $Q\subset M$ is \emph{periodic} if there is $k\in {\mathbb N}$ such that $S^k(x)=x$ for any $x\in Q$ and all iterations $S^l|_Q$ are homeomorphisms for each $1\le l\le k$.
\end{definition}

For piecewise translation maps (but not for piecewise isometries) this implies that all points of $Q$ are periodic.

\begin{definition} We say that a subset $Q\subset M$ is \emph{eventually periodic} if there is $m\in {\mathbb N}$ such that $S^m(Q)$ is a periodic subset.
\end{definition}
	
\begin{lemma}
Let $S$ be a piecewise translation map. If $x_0\in M{\setminus}\Omega$ is a periodic point of $S$, there exists a periodic ball that contains the point $x_0$.
\end{lemma}

\noindent\textbf{Proof.} Let $k \ge 1$ be such that $S^k(x_0)=x_0$. Since $x_0\notin M{\setminus}\Omega$ there exists a $\delta> 0$ such that if $B_\delta$ is the open ball of radius $\delta$ centered at $x_0$, then $S^l(B_\delta)$ is the open ball of radius $\delta$ centered at $S^l(x_0)$ for each $0\le l\le k$. In particular, $S^k(B_\delta)=B_\delta$, showing that the ball $B_\delta$ is periodic. $\square$

A similar statement is true for eventually periodic points.

\begin{definition} We call a domain $J\subset M$ \emph{tough} if all maps $S^k$, $k\in {\mathbb N}$, are continuous on $J$. In one-dimensional case we use the classical notion \emph{homterval} for tough intervals.
\end{definition}

Evidently, this means that discontinuity points do not belong to $S^k(J)$ for all $k\in {\mathbb N}$. We call a tough domain \emph{maximal} if it is not a proper subset of another tough domain. Evidently if $J_0$ is a tough domain then all sets $S^k(J_0)$ ($k\ge 0$) are tough domains.

\begin{lemma} For any interval translation map with rational values $c_j$, all points of ${\mathbb T}^1{\setminus}\Omega$ are eventually periodic and the set ${\mathbb T}^1{\setminus}\Omega$ is a finite union of homtervals.
\end{lemma}

Now we go back to the general case of piecewise isometries.

\begin{proposition} The following statements hold.
\begin{enumerate}
\item Let $J_\alpha, \alpha\in A$ be a family of tough domains such that $$\bigcap_{\alpha \in A} J_\alpha\neq \emptyset.$$
Then 
$$J =\bigcup_{\alpha\in A} J_\alpha$$
is a tough domain. 
\item If $J_1$ and $J_2$ are tough domains such that $J_1\bigcap J_2\neq \emptyset$ then $J_1 \bigcap J_2$ is a tough domain. 
\end{enumerate}
\end{proposition}

\begin{lemma} Any tough domain $J_0$ is a subset of the uniquely defined maximal tough domain.
\end{lemma}

\noindent\textbf{Proof.} This domain may be defined as the union of all tough domains $U\supset J_0$. $\square$

\begin{lemma} Any tough domain is eventually periodic.
\end{lemma}

\noindent\textbf{Proof.} Consider a tough domain $J_0$. Since the Lebesgue measure of the manifold $M$ is finite and $\mathrm{Leb}(S^k(J_0))=\mathrm{Leb}(J_0)>0$ for all positive $k$, there exist $k,l\in {\mathbb N}$, $k> l$ such that $S^l(J_0)\bigcap S^k(J_0)\neq \emptyset$. Take $D$ -- the maximal tough domain that contains $S^l(J_0)$. Then $S^{k-l}(D) \bigcup D$ is also a tough domain which means that $S^{k-l}(D)=D$. Since $S^{k-l}|D$ is a piecewise isometry, we are done. $\square$ 

\begin{definition} We call a map $S$ \emph{generic} if it does not have any periodic $($tough$)$ domains. 
\end{definition}

Given a partition $\{[t_j,t_{j+1})\}$ of the circle ${\mathbb T}^1$, we consider the space of all piecewise translation maps, varying values of shifts $c_j$. Observe that all translation maps form a compact subset of a Euclidean space. 

\begin{corollary} For the interval translation map of Example 1.2 given a fixed number of domains, for a generic map $S$, the set $\Omega$ is dense in ${\mathbb T}^1$.
\end{corollary}

\noindent\textbf{Proof.} If the set $\Omega$ is not dense, there exists a nontrivial homterval. So, some of shift vectors must be rationally dependent. The set of all parameters $t_j$ and $c_j$ with rationally independent $c_j$ is generic.
$\square$
 
\section{Invariant measures for interval translation maps.}

In this section, we consider interval translation maps only. Particularly, we set $M={\mathbb T}^1$ everywhere in this section. 

Later on we consider Borel probability invariant measures only. We say that a measure $\mu$ is invariant with respect to the map $S$ if $\mu(S^{-1}(A))=\mu (A)$ for any measurable set $A$. Since the map $S$ is discontinuous, we cannot appeal to the Krylov-Bogolyubov theorem for existence of invariant measures. We cannot say that the Lebesgue measure is invariant. However, as we demonstrate below, invariant measures exist for any map $S$.

\begin{definition} A measure is called \emph{non-atomic} if the measure of every singleton is zero. 
\end{definition}

Observe that if a map does not have any periodic points, every invariant probability measure is non-atomic. Supports of non-atomic measures are uncountable.

\begin{theorem}\label{th1d} Any interval translation map admits Borel probability non-atomic invariant measures.
\end{theorem}

\noindent\textbf{Proof.} Let the map $S$ have periodic domains and let $J_0$ be one of them. Take the minimal positive value $m$ such that $S^m(J_0)=J_0$. Then the invariant measure $\mu^*$ can be constructed as a renormalised restriction of the Lebesgue measure to the set $J_0\bigcup S(J_0) \bigcup \ldots \bigcup S^{m-1}(J_0)$.

So, we may assume that the map $S$ is generic. Consider sequences $\{t^m_k\}$ and $\{c^m_k\}$ of rational numbers that converge to $t_k$ and $c_k$ respectively. Let $S_m$ be corresponding mappings.

\begin{lemma}\label{elA} The approximating sequences $\{t^m_k\}$ and $\{c^m_k\}$ may be selected so that the following statement is true. Assume that numbers $i,j\in \{1,\ldots,n\}$: $i\neq j$, and 
$l_1,\ldots, l_n \in {\mathbb Z}^+$ are such that 
\begin{equation}\label{linear}
t_j-t_i=\sum_{k=1}^n l_k c_k.
\end{equation}
Then for any $m\in {\mathbb N}$ 
\begin{equation}\label{linear1}
t^m_j-t^m_i=\sum_{k=1}^n l_k c^m_k.
\end{equation}
\end{lemma}

\noindent\textbf{Proof.}  Let $L$ be the space of all values $t_j$ and $c_j$ satisfying all equations \eqref{linear} that are true for the considered values of parameters.  Since all coefficients $l_k$ are integers, points with all rational coordinates are dense in $L$. 
$\square$

\begin{lemma}\label{elB}
For any $r\in {\mathbb N}$, $i,j=0,\ldots, n$ such that $$S^r (t_i+0)=t_j \quad \mbox{or} \quad S^r(t_i-0)=t_j,$$ 
there exists $m_r\in {\mathbb N}$ such that $S^r_m(t^m_i+0)= t_j^m$ or, respectively, $S^r_m(t_i^m+0)= t_j^m$\, for any $m\ge m_r$. 
\end{lemma}

\noindent\textbf{Proof.} We select approximating sequences $t^m_k$ and $c_k^m$ so that given $i, j, l_1,\ldots, l_n$ all equalities \eqref{linear1} are satisfied for all $m$ provided the corresponding equality \eqref{linear} is true. Take $m_1$ so big that  
for any $j\in \{1,\ldots,n\}$ the inclusion
$$
S(t_j+0) \in (t_i,t_{i+1})
$$
implies $S_m(t^m_j+0) \in (t^m_i,t^m_{i+1})$ and the equality 
$$S(t_j+0) = t_i$$
implies $S_m(t^m_j+0) =t_i$ for all $m\ge m_1$. 

For such values of $m$, the boundaries of continuity segments for the map $S_m^2$ converge to ones of the map $S^2$. All corresponding shifts (that form a subset of $\{c^m_i+c^m_j \mod 1: i,j=1,\ldots,n\}$) converge to corresponding shifts of the map $S^2$. So, there exists $m_2$ such that for any $j\in \{1,\ldots,n\}$ the inclusion
\begin{equation}\label{insegment}
S^2(t_j+0) \in (t_i,t_{i+1})
\end{equation}
implies $S^2_m(t^m_j+0) \in (t^m_i,t^m_{i+1})$
and the equality 
$$S^2(t_j+0) = t_i$$
implies $S^2_m(t^m_j+0) =t_i$ for all $m\ge m_2$ \ldots

To finish the proof it suffices to repeat the similar procedure $r$ times. $\square$

Let $$\Xi_m=\bigcap_{k=1}^\infty S_m^k({\mathbb T}^1).$$ 
These sets have positive Lebesgue measures, since each of them is a union of a finite number of arcs. We introduce measures $\mu_m$ as renormalisations of the Lebesgue measure, restricted to  $\Xi_m$. These measures are invariant w.r.t. mappings $S_m$.

\begin{lemma}\label{l_tk} Let the map $S$ be generic. For any $\varepsilon>0$ there exists $\delta>0$ such that 
\begin{equation}\label{smallmeasure}
\limsup\limits_{m\to \infty} \mu_m(t_k-\delta,t_k+\delta)<\varepsilon
\end{equation}
for any $k=1,\ldots,n$.
\end{lemma}

\noindent\textbf{Proof.} Let the statement of the lemma be wrong. Assume, without loss of generality, that there exist a number $\varepsilon>0$, a fixed number $k$ and sequences 
$\{m_r\}\subset {\mathbb N}$ and $\delta_r \to 0$ such that 
$$\mu_{m_r} (t_k-\delta_r,t_k+\delta_r)\ge\varepsilon.$$
Let $N\in {\mathbb N}$ be such that $1/N<\varepsilon$. Then for any $r$ there exist numbers $0\le i<j \le N$ such that 
$$S_{m_r}^{-i} (t_k-\delta_r,t_k+\delta_r) \bigcap S_{m_r}^{-j} (t_k-\delta_r,t_k+\delta_r) \neq \emptyset.$$
Let $l=j-i$. Then for any $r\in {\mathbb N}$ there exists a point $x_r\in (t_k-\delta_r,t_k+\delta_r)$ such that $S_{m_r}^l(x_r) \in (t_k-\delta_r,t_k+\delta_r)$. Since we could, without loss of generality, take the same number $l$ for all values of $r$, at least one of two statements is true: either $S^l(t_k-0)=t_k$ or $S^l(t_k+0)=t_k$. Both cases imply existence of a periodic domain that contradicts assumptions of the lemma. $\square$

\begin{corollary}\label{cor} Assume that there exists a weak-$*$ limit $\mu^*$ of measures $\mu_m$. Then
\begin{enumerate}
\item $\mu^*(\{t_k\})=0$ for any $k=1,\ldots,n$.
\item for any $k=1,\ldots,n$ and any continuous function $\varphi:[t_{k-1},t_k]\to {\mathbb R}$ 
$$\lim_{m\to\infty} \int_{[t_{k-1},t_k]} \varphi\, d\mu_m= \int_{[t_{k-1},t_k]} \varphi\, d\mu^*.$$
\end{enumerate}
\end{corollary}

\noindent\textbf{Proof.} The first statement is evident. To prove the second one we assume without loss of generality that 
$$\max_{[t_{k-1},t_k]} |\varphi(t)|=1.$$ Then given an $\varepsilon>0$ we take a $\delta>0$ so that $\mu_m([t_{k-1},t_{k-1}+\delta]\bigcup [t_k-\delta,t_k])<\varepsilon$ for any $m$ and 
$\mu^*([t_{k-1},t_{k-1}+\delta]\bigcup [t_k-\delta,t_k])<\varepsilon$ (the last inequality may be satisfied due to the first statement of the lemma).

Take a non-negative continuous function $\eta: {\mathbb T}^1 \to [0,1]$ such that $\eta (t)=0$ for any $t \le t_{k-1}$ and for any $t \ge t_k$ and $\eta(t)=1$ for any $t\in [t_{k-1}+\delta, t_k-\delta]$. Let $\psi(t)=\eta(t)\varphi(t)$ if $t\in [t_{k-1},t_k]$ and $\psi(t)=0$ otherwise. This function is continuous. Then $$\lim_{m\to\infty} \int_{[t_{k-1},t_k]} \psi\, d\mu_m= \int_{[t_{k-1},t_k]} \psi\, d\mu^*.$$
This means that 
$$\limsup_{m\to\infty} \left| \int_{[t_{k-1},t_k]} \varphi\, d\mu_m- \int_{[t_{k-1},t_k]} \varphi\, d\mu^*\right|\le 2\varepsilon.$$
Since the value $\varepsilon$ can be taken arbitrarily small, this proves the statement of the corollary. $\square$

\begin{lemma}\label{l_tech} For any continuous function $\varphi : {\mathbb T}^1\to {\mathbb R}$, we have 
$$\int_{{\mathbb T}^1} (\varphi(S_m(t)) -\varphi(S(t))) \, d\mu_m(t) \to 0$$
as $m\to\infty$.
\end{lemma}

\noindent\textbf{Proof.} Fix a positive $\varepsilon>0$. For this $\varepsilon$, select $\delta>0$ so that all inequalities \eqref{smallmeasure} are satisfied.
We can find $m_0\in {\mathbb N}$  such that 
$$|\varphi(S_m(t))-\varphi(S(t))|<\varepsilon$$ 
for any $m\ge m_0$ and
$$t\notin \bigcup_{k=1}^n [t_k-\delta,t_k+\delta]$$ 
(here we take $m_0$ so big that $|t_k-t_k^m|<\delta/2$ for all $m\ge m_0$, $k=1,\ldots,n$).
So, 
$$\int_{{\mathbb T}^1{\setminus} \bigcup [t_k-\delta,t_k+\delta]} |\varphi(S_m(t))-\varphi(S(t))|\, d\mu_m(t) \le \varepsilon.
$$ 
For any $k$, we have 
$$\begin{array}{c}
\int_{\bigcup [t_k-\delta,t_k+\delta]} |\varphi(S_m(t))-\varphi(S(t))|\, dt \le \\
2 \max_{t\in [t_k-\delta,t_k+\delta]} |\varphi (t)| \mu_m\left(\bigcup [t_k-\delta,t_k+\delta]\right)\le 2 \max |\varphi (t)|\varepsilon.
\end{array}
$$
Since the value $\varepsilon$ can be taken arbitrarily small, this finishes the proof of the lemma. $\square$ 

Now, we continue the proof of the theorem. 

The set of Borel probability measures in ${\mathbb T}^1$ is compact in the weak-$*$ topology. Without loss of generality, we may assume that the sequence of invariant measures $\mu_k$ weakly-$*$ converges to a measure $\mu^*$. 

Let us prove that $\mu^*$ is invariant w.r.t. $S$. 

Fix a continuous function $\varphi:{\mathbb T}^1\to {\mathbb R}$. We need to prove that 
$$\int_{{\mathbb T}^1} \varphi\, d\mu^*=\int_{{\mathbb T}^1} \varphi \circ S\, d\mu^*.$$

We have
\begin{equation}\label{9.4}
\begin{array}{c}
\left|\int\limits_{{\mathbb T}^1} (\varphi- \varphi \circ S)\, d\mu^*\right|\le 
\left|\int\limits_{{\mathbb T}^1} (\varphi- \varphi \circ S)\, d\mu^*-\int\limits_{{\mathbb T}^1} (\varphi- \varphi \circ S)\, d\mu_m\right|+ \\[10pt]
\left|\int\limits_{{\mathbb T}^1} (\varphi- \varphi \circ S_m)\, d\mu_m\right|+\left|\int\limits_{{\mathbb T}^1}( \varphi\circ S_m - \varphi \circ S)\, d\mu_m\right|.
\end{array}
\end{equation}
The first term in the right hand side of Eq.\, \eqref{9.4} tends to $0$ due to Corollary \ref{cor}, the second one is zero since measures $\mu_m$ are $S_m$-invariant and the third one tends to zero by Lemma  \ref{l_tech}. So, the right hand side of Eq.\, \eqref{9.4} is zero. 

The obtained measure is non-atomic since the generic map does not have any periodic points out of $\Omega$ (see Lemma 2.1) and $\mu^*(\Omega)=0$. 
Recall that the case of maps with periodic points was studied in the beginning of the proof. $\square$

Let $\Xi=\overline {\bigcap_{k=1}^\infty S^k({\mathbb T}^1)}$.

\begin{corollary} For any map $S$ the set $\Xi$ is uncountable. 
\end{corollary}

\noindent\textbf{Proof}. The support of the non-atomic invariant measure that exists by Theorem 3.1 is a subset of $\Xi$. $\square$

\begin{lemma} Let $\mu^*$ be the non-atomic invariant measure for an interval translation map $S$ that exists by Theorem \ref{th1d}. Then
\begin{equation}\label{musa}
\mu^*(S([a,b]))=\mu^*([a,b]) 
\end{equation}
for any segment $[a,b]\subset {\mathbb T}^1\, {\setminus}\, \{t_1,\ldots,t_n\} $.
\end{lemma}

\noindent\textbf{Proof.}  Similarly to the proof of Theorem \ref{th1d}, we approximate the map $S$ by maps $S_m$ for which all points are eventually periodic and weakly-$*$ approximate the measure $\mu$ by continuous measures $\mu_m$.  All maps $S_m$ are  invertible almost everywhere on supports of their invariant measures. So, 
\begin{equation}\label{mumsma}
\mu_m(S_m([a,b]))=\mu_m ([a,b])
\end{equation}
for any $m\in {\mathbb N}$ and any segment $[a,b]$ (in fact, the similar statement is true for any measurable set). Let us prove that we can take the limit in \eqref{mumsma}. Given $\varepsilon>0$, we consider $\delta>0$ such that 
$$\mu_m ([a-\delta,a]\bigcup [b,b+\delta])<\varepsilon$$ and 
$\mu^* ([a-\delta,a]\bigcup [b,b+\delta])<\varepsilon$. 

Consider a continuous function $\eta: {\mathbb T}^1 \to [0,1]$ that equals 1 at points of $[a,b]$ and equals 0 out of $[a-\delta,b+\delta].$
Then, 
$$\int_{{\mathbb T}^1} \eta\, d\mu_k \to \int_{{\mathbb T}^1} \eta\, d\mu^*$$
and, consequently,
$$\liminf \mu_m ([a,b])-\varepsilon \le \mu^*([a,b]) \le \limsup \mu_m ([a,b])+\varepsilon$$
and, therefore 
$$\mu^*([a,b])=\lim \mu_m([a,b]).$$

Observe that all maps $S_m$ are continuous on $[a,b]$ for big values of $m$. Then there exist constants $c_m$ and $c$ such that $S_m([a,b])=[a+c_m,b+c_m]$ and $S([a,b])=[a+c,b+c]$. Then, similarly to the previous step, we can prove that  
$$\lim \mu_m(S_m([a,b])) = \mu^*(S([a,b]).$$
$\square$

Taking limits in \eqref{musa}, we could prove the similar statement for segments $[a,b]\subset [t_k,t_{k+1}]$ even in case when $a=t_k$ or/and $b=t_{k+1}$. Here we use the fact that the measure $\mu^*$ is non-atomic.

We show that any interval translation map endowed with a non-atomic invariant measure is metrically equivalent to an interval exchange map of the segment $[0,1]$.

\begin{theorem} Let $\mu^*$ be the non-atomic invariant measure for an interval translation map $S$ that exists by Theorem \ref{th1d}. Then the restriction $S|_{{\mathrm{supp}\,} \mu^*}$ is metrically equivalent to an interval exchange map $T:[0,1]\to [0,1]$ with the Lebesgue measure. The semi-conjugacy map is one-to-one everywhere, except a countable set.
\end{theorem}

\noindent\textbf{Proof.}  Consider the function $h(x)=\mu^*([0,x])$. This function restricted to ${\mathrm{supp}\,} \mu^*$  is monotonous. The measure $\mu^*$ is non-atomic, so $h$ is continuous. The equality $h(x)=h(y)$ implies $\mu^* ([x,y])=0$ so there is a countable set 
$D\subset {\mathrm{supp}\,} \mu^*$ such that $h|_{{\mathrm{supp}\,} \mu^* {\setminus} D}$ is injective and $\# h^{-1}(y)=2$ for any $y\in h(D)$. Now we define the map $\hbar(x):=\max h^{-1}(x):[0,1]\to {\mathrm{supp}\,} \mu^*$ that is a right inverse to $h$ and set
\begin{equation}\label{deft}
T(x):=h (S(\hbar(x)).
\end{equation}
Then, by definition we have the semi-conjugacy: $T\circ h=h\circ S$ out of the countable set $D$.

Now let $0= t_0<\ldots < t_n=1$ be the points such that the maps $S|_{[t_j,t_{j+1})}$ are translations: $S(x)=x+c_j$  for all $j=0,\ldots, n$. For any $j$ we denote $\tau_j=\mu^*([0,t_j])$, so 
$\tau_j=h(t_j)$ if we naturally extend the map $h$ to $[0,1]$.

Take two points $x<y$ in a segment $I_j:=[\tau_j,\tau_{j+1})$ if the segment is non-empty. By \eqref{deft} we have
$$
T(y)-T(x)=h(S(\hbar (y))-h(S(\hbar(x)))
$$
The map $\hbar$ is monotonous and $\hbar(x),\hbar (y)$ belong to the same segment $[t_j,t_{j+1})$ where $S$ is the shift. So, using Eq. \eqref{musa}, we have
$$T(y)-T(x)=\mu^* [S(\hbar(x)),S(\hbar(y))]=\mu^* [\hbar(x),\hbar(y)]= h(\hbar(y))-h(\hbar(x))=y-x.$$

So, $T$ is an interval translation map. To finish the proof, it suffices to demonstrate that $T$ preserves the Lebesgue measure. Observe that ${\mathrm{Leb}\,} (h(I))=\mu^*(I)$ for any segment $I$. A similar statement is true for any set that is a finite union of segments. We have
$$\begin{array}{c} 
{\mathrm{Leb}\,} T^{-1} ([x,y))= {\mathrm{Leb}\,}(\hbar^{-1}\circ S^{-1}\circ h^{-1}([x,y)))={\mathrm{Leb}\,} (h(S^{-1}\circ h^{-1}([x,y))))=\\
\mu^*(S^{-1}\circ h^{-1}([x,y)))=\mu^* (h^{-1}([x,y)))=\mu^* [\hbar(x),\hbar(y))=y-x.
\end{array}$$
Then ${\mathrm{Leb}\,}(T^{-1}(A))={\mathrm{Leb}\,} (A)$ for any measurable set $A$. 

If the map $T$ is not an interval exchange map, there must be two segments whose images coincide. This contradicts to the invariance of the Lebesgue measure. $\square$

Let us recall some standard definitions from Topological Dynamics.

\begin{definition} We call a point $x$ \emph{recurrent (Poisson stable)} with respect to the map $S$ if there exists an increasing sequence $\{m_k\in {\mathbb N}\}$ such that $S^{m_k}(x) \to x$.
\end{definition}

\begin{definition} We call a point $x$ \emph{nonwandering} with respect to the map $S$ if for any neighbourhood $U$ of the point $x$ there is a point $p\in U$ and a number $k\in {\mathbb N}$ such that $p,S^k(p)\in U$. 
\end{definition}

Observe that the last definition works even for points of the discontinuity set where the map may be undefined.

Recall that the closure of the set of all recurrent points is a subset of the set $\Omega(S)$ of all non-wandering points.The next statement is well-known for continuous maps, we give a proof for interval translation maps.

\begin{theorem}\label{thpois} Let $\mu^*$ an the invariant measure for the mapping $S$. The recurrent points of $S$ are dense in ${\mathrm{supp}}\, \mu^*$.
\end{theorem}

\noindent\textbf{Proof.} Let $\{U_j\}_{j\in {\mathbb N}}$ be a countable base of the topology in ${\mathrm{supp}}\, \mu^*$. For any $j$ consider the set 
$$N_j:=\{x\in U_j: \# \{m\in {\mathbb N}: S^m(x)\in U_j\}<\infty\}.$$
By Poincar\'e Recurrence Theorem, $\mu(N_j)=0$ for any $j$, so $\mu (R)=1$ where 
$$R={\mathrm{supp}}\, \mu^* {\setminus} \bigcup_{j=1}^\infty N_j.$$

For any $x\in R$ and any $\varepsilon>0$ there exists $m\in {\mathbb N}$ such that $|S^m(x)-x|<\varepsilon$. Then $x$ is a recurrent point. The set $R$ is dense in ${\mathrm{supp}}\, \mu^*$ by definition of the support of the measure. $\square$ 

So, the set of recurrent points must be uncountable for interval translation maps. 

\begin{remark} In the proof of Theorem \ref{th1d} it was very important that we consider a one-dimensional map. Existence of an invariant measures for all piecewise isometries (or, even, for all piecewise translation maps) is an open question. The only non-trivial result in this area, author knows, follows from the main theorem of the preprint \cite{T}. There, a piecewise translation map $S:M\to M$ is considered where $M=\bigcup_{j=1}^n M_j$ is a subset of ${\mathbb R}^d$. Let cardinality $n$ of the partition $\{M_j\}$ be such that $n=d+1$ and $c_j$ be the translation vectors. It was proved that the map $S$ is finite provided vectors $c_1,\ldots c_{d+1}$ be rationally independent. In this case the Lebesgue measure of the set 
$$\Omega=\bigcap_{m=1}^\infty S^m (M)$$
is non-zero. So, the restriction of the Lebesgue measure to $\Omega$ is a non-atomic invariant measure. It is an easy exercise that an invariant measure exists if $n<d+1$ and the vectors $c_1, \ldots, c_n$ are rationally independent (we leave the proof to readers).
\end{remark}

\section{Invariant measures for Piecewise Continuous Maps} 

Here we formulate and prove a result on existence of invariant measures for piecewise continuous maps. Once again, we consider a $d$-dimensional compact Riemannian manifold $(M,{\mathrm{dist}\,})$ or a bounded set $M\subset {\mathbb R}^d$ that can be represented as a finite union of closures of pairwise disjoint domains $M_j$ ($j=1,\ldots, n$). Consider a map $T:M\to M$ where restrictions $T|_{M_j}$ are continuous maps.

Let $H$ be the (compact) discontinuity set of $T$ and $H_\varepsilon=\{x\in M: {\mathrm{dist}\,} (x,H)<\varepsilon\}$. Let $\Omega:=\bigcup_{j=0}^\infty T^{-j} (H)$.

\begin{theorem} Assume that there exists a point $x\notin \Omega$ such that  
\begin{equation}\label{eliminf}
\lim_{\varepsilon\to 0} \limsup_{m\in {\mathbb N}} \dfrac{1}{m}\# \{1 \le k\le m: T^k(x) \in H_\varepsilon\}=0.
\end{equation}
Then the map $T$ admits a Borel probability invariant measure.
\end{theorem}

\noindent\textbf{Proof.} We use the idea of the proof of the Krylov-Bogolyubov theorem. We fix the point $x\in M{\setminus} \Omega$ and consider the sequence of measures 
$$\mu_m:=\dfrac{1}{m} \sum_{i=0}^{m-1} \delta(T^i(x))$$
where $\delta(y)$ is the Dirac measure at the point $y$. All these measures are Borel probabilities on $M$. The space of such measures is compact in the weak-$*$ topology hence we could take a subsequence $\mu_{m_k}$, weakly-$*$ converging to a measure 
$\mu^*$. The assumption \eqref{eliminf} guarantees that given $\sigma>0$ there exist $m_0\in {\mathbb N}$ and $\varepsilon>0$ such that $\mu_m(H_\varepsilon)\to 0$ for any $m\ge m_0$.
Therefore $\mu^*(H)=\mu^*(\Omega)=0$. Let $T^\#$ be the pushforward operator on the space of Borel probability measures: $T^\#(\mu)=\nu$ if $\nu(A)=\mu(T^{-1}(A))$ for any measurable set $A$. Observe that for any $m$
$$T^\# \mu_m-\mu_m=\dfrac1{m} (\delta (T^m(x))-\delta(x)).$$
So, 
\begin{equation}\label{elim}
\lim\limits_{k\to\infty} T^\#(\mu_{m_k})=\mu^*. 
\end{equation}

\begin{lemma} $\lim\limits_{k\to\infty} T^\#(\mu_{m_k})=T^\# (\mu^*)$.
\end{lemma}

\noindent\textbf{Proof.} Fix a continuous function $\varphi: M\to {\mathbb R}$. We need to prove that 
\begin{equation}\label{last}
\lim \int_M \varphi\circ T \, d\mu_{m_k} = \int_M \varphi \circ T \, d\mu^*.
\end{equation}

Fix an $\varepsilon>0$ and consider a continuous function $\eta_\varepsilon:M\to [0,1]$ that equals $0$ in $H_\varepsilon$ and equals $1$ out of $H_{2\varepsilon}$. Functions $\varphi(T(x))\eta_\varepsilon(x)$ are continuous for any $\varepsilon>0$, so
$$\lim \int_M \varphi\circ T\cdot \eta_\varepsilon \, d\mu_{m_k} = \int_M \varphi \circ T \cdot \eta_\varepsilon\, d\mu^*.$$

Given a $\kappa>0$ and a function $\varphi$, we could select $m_0$ and $\varepsilon$ so that 
$$\left| \int_M \varphi\circ T\cdot \eta_\varepsilon \, d\mu_{m_k} - \int_M \varphi\circ T\, d\mu_{m_k}\right| \le \kappa$$
for any $m\ge m_0$ and 
$$\left| \int_M \varphi\circ T\cdot \eta_\varepsilon \, d\mu^* - \int_M \varphi\circ T\, d\mu^*\right| \le \kappa.$$
This proves the lemma. Then statement of the lemma and formula \eqref{elim} imply that the measure $\mu^*$ is invariant.  $\square$
 
\begin{corollary} If the discontinuity set $H$ of a map $T$ does not contain nonwandering points, the map $T$ admits a Borel probability invariant measure.
\end{corollary} 

The last result of the paper concerns numerical approximations of invariant measures for discontinuous maps. Also, it gives a weak version of semicontinuity of non-atomic invariant measures with respect to parameters of ITMs. 

\begin{theorem} Let $S: {\mathbb T}^1 \to {\mathbb T}^1$ be a piecewise continuous map with a finite discontinuity set $H$. Let a sequence of interval translation maps $S_k$ converge to $S$ uniformly on compact subsets of ${\mathbb T}^1\,{\setminus}\, H$. 
Assume that a sequence of $S_k$-invariant Borel probability measures $\mu_k$ converges weakly-$*$ to a measure $\mu^*$ and
\begin{equation}\label{ecrt} \mu^*(H)=0.
\end{equation} 
Then the measure $\mu^*$ is $S$-invariant.
\end{theorem}

\noindent\textbf{Proof}. Take a continuous function $\varphi:{\mathbb T}^1 \to {\mathbb R}$ and $\varepsilon>0$. Given a $\delta>0$, let  $\eta_\delta \in C^0({\mathbb T}^1\to [0,1])$ be a continuous function such that $\eta_\delta=0$ in $U_\delta(H)$ -- the $\delta$-neighbourhood of $H$ and $\eta_\delta(x)=1$ out of the $2\delta$-neighbourhood of $H$. 

Observe that given $\kappa>0$, we can take $\delta>0$ so that $\limsup_k \mu_k (U_\delta(H)) <\kappa$. To see this, it suffices to consider convergence 
$$\int_{{\mathbb T}^1} (1-\eta_\delta) \, d\mu_k\to \int_{{\mathbb T}^1} (1-\eta_\delta) \, d\mu^*$$
and take into account the fact that $\mu^*(U_{2\delta}(H)) \to 0$ as $\delta\to 0$. 

Since the measures $\mu_k$ are $S_k$-invariant, we have 
$$\int_{{\mathbb T}^1} \varphi \, d\mu_k= \int_{{\mathbb T}^1} \varphi \circ S_k\, d\mu_k$$
for all $k\in {\mathbb N}$. Given $\varepsilon>0$, we take $\delta>0$ so small that  
$$\max_{x\in {\mathbb T}^1}|\varphi (x)| \limsup_{k\to \infty} \mu_k (U_{2\delta}(H)) <\varepsilon$$
(we can do this since the function $\varphi$ is continuous). Then evidently
$$\max_{x\in {\mathbb T}^1}|\varphi (x)| \limsup_{k\to \infty} \mu^* (U_{2\delta}(H)) <\varepsilon.$$

Then $$\limsup_k \left| \int_{{\mathbb T}^1} \varphi \cdot \eta_\delta \, d\mu_k- \int_{{\mathbb T}^1} \varphi \circ S_k \cdot \eta_\delta\, d\mu_k\right|\le 2\varepsilon.$$

Since $S_k$ converges to $S$ uniformly on $ {\mathbb T}^1{\setminus} U_\delta(H)$, the latter formula implies 
$$
\limsup_k \left| \int_{{\mathbb T}^1} \varphi \cdot \eta_\delta \, d\mu_k- \int_{{\mathbb T}^1} \varphi \circ S \cdot \eta_\delta\, d\mu_k\right|\le 2\varepsilon.
$$
Both functions $\varphi \cdot \eta_\delta$ and $\varphi \circ S \cdot \eta_\delta$ are continuous of ${\mathbb T}^1$, so, we may proceed to limit in \eqref{last}. So,
$$\left| \int_{{\mathbb T}^1} \varphi \cdot \eta_\delta \, d\mu^*- \int_{{\mathbb T}^1} \varphi \circ S \cdot \eta_\delta\, d\mu^* \right|\le 2\varepsilon$$
and, consequently 
$$\left| \int_{{\mathbb T}^1} \varphi \, d\mu^*- \int_{{\mathbb T}^1} \varphi \circ S \, d\mu^* \right|\le 4\varepsilon.$$
Since the continuous function $\varphi$ is arbitrary and $\varepsilon>0$ is arbitrarily small, the measure $\mu^*$ is $S$-invariant. $\square$ 

\section{Conclusion} Here we list once again the principal results of the paper.  
\begin{enumerate} 
\item Any interval translation map admits a Borel non-atomic probability invariant measure. 
\item Interval translation maps endowed with this measure are metrically equivalent to interval exchange maps with the Lebesgue measure.
\item Any piecewise continuous map without nonwandering points on the discontinuity set admits a Borel probability invariant measure.
\item If a sequence of interval translation maps approximates a piecewise continuous map and the sequence of corresponding invariant measures, weakly converges to a measure, the latter measure is invariant with respect to the piecewise continuous map.
\end{enumerate}

\section*{Acknowledgements} 
The work was partially supported by RFBR grant 18-01-00230-a. The author thanks M.\,Arnold, M.\,Artigiani, V.\,Avrutin, N.\,Begun, D.Rachinsky, D.\,Volk and many others for their support and precious advices. Also, he is grateful to the anonymous referees for their remarks. This work is dedicated to author's wife Maria and daughter Elizaveta, who were inspiring him all this time.

\end{document}